\documentclass[11pt]{article}

\usepackage{amssymb,amsmath,amsthm,pslatex,graphicx}
\usepackage{mathrsfs}
\usepackage{epsfig}
\usepackage{tikz}
\usepackage{enumerate}

\swapnumbers
\newtheorem{theorem}{Theorem}[section]
\newtheorem{lemma}[theorem]{Lemma}

\renewcommand\proof{\noindent\textsl{Proof. }}
\newcommand\sqr[2]{{\vbox{\hrule height.#2pt
    \hbox{\vrule width.#2pt height#1pt \kern#1pt
        \vrule width.#2pt}\hrule height.#2pt}}}

\renewcommand\qed{%
\ifmmode\eqno\sqr53
\else\nolinebreak\ \hfill\sqr53\medbreak\fi}

\begin{document}

\title{Sperner partition systems}

\author{
    P.~C.  Li\thanks{
    Research supported by NSERC Discovery Grant 250389-06
   } \\
    Dept. of Computer Science\\
    University of Manitoba\\
    Winnipeg, Manitoba\\
    Canada R3T 2N2
\and
    Karen Meagher\thanks{
    Research supported by NSERC Discovery Grant 341214-08
   } \\
    Dept. of Mathematics and Statistics  \\ 
    University of Regina \\
    Regina, Saskatchewan \\
    Canada S4S 0A2   
 }

\maketitle
\begin{abstract} 
  A \textsl{Sperner $k$-partition system} on a set $X$ is a set of
  partitions of $X$ into $k$ classes such that the classes of the
  partitions form a Sperner set system (so no class from a partition
  is a subset of a class from another partition). These systems were
  defined by Meagher, Moura and Stevens in~\cite{MMS} who showed that
  if $|X| = k \ell$, then the largest Sperner $k$-partition system has
  size $\frac{1}{k}\binom{|X|}{\ell}$. In this paper we find bounds on
  the size of the largest Sperner $k$-partition system where $k$ does
  not divide the size of $X$, specifically, we give an exact bound
  when $k=2$ and upper and lower bounds when $|X| = 2k+1$, $|X|=2k+2$
  and $|X| = 3k-1$.
\end{abstract}

\maketitle

\section{Introduction}

This paper continues the work from~\cite{MMS} that established Sperner
type theorems for partitions. Two subsets from a set $X$ are call
\textsl{incomparable} if neither set is contained in the other. A
\textsl{Sperner set system} on a set $X$ is a collection of subsets of
$X$ such that any two are incomparable. In 1928 Sperner gave the exact
size of the largest Sperner set system.

\begin{theorem}(Sperner's Theorem~\cite{sperner:28}) A Sperner set system on a
  set $X$ with cardinality $n$ has at most $\binom{n}{\lfloor n/2
    \rfloor}$ sets.  The only Sperner set systems that meet this bound
  are the collection of all sets of size $\lfloor n/2 \rfloor$ or the
  collection of all sets of size $\lceil n/2 \rceil$.
\end{theorem}

A partition of a set $X$ is a disjoint collection of non-empty
subsets, called classes, of $X$ whose union is $X$. A partition with
exactly $k$ classes is called a $k$-partition. A \textsl{Sperner
  $k$-partition system} on a set $X$ is a set of $k$-partitions of $X$
such that the classes of the partitions form a Sperner set
system. This property is equivalent to requiring that no class from
any partition is a subset of a class of another partition. A Sperner
$k$-partition system can be considered to be a \textsl{resolvable Sperner
  set system}.  Figure~\ref{fig1} shows an example of a Sperner
$3$-partition system on $n=7$ points.

\begin{figure} 
\label{fig1}
\begin{center}
\begin{tabular}{c|l} 
partition & classes \\ \hline
$1$ & $\{0,1\},\;\{2,3\},\;\{4,5,6\}$ \\ 
$2$ & $\{0,2,6\},\;\{1,5\},\;\{3,4\}$ \\
$3$ & $\{0,3\},\;\{1,4\},\;\{2,5,6\}$ \\
$4$ & $\{0,4\},\;\{1,2,6\},\;\{3,5\}$ \\ 
$5$ & $\{0,5\},\;\{1,3\},\;\{2,4,6\}$
\end{tabular}
\end{center}
\caption{A Sperner $3$-partition system on $\{0,...,6\}$ with $5$
  partitions.}
\end{figure}

These systems were introduced by Meagher, Moura, Stevens~\cite{MMS}
who established bounds on the number of partitions in such a system.
Given $n,k$, let $SP(n,k)$ denote the size of the largest (in terms of
number of partitions) Sperner $k$-partition system. In~\cite{MMS} the
exact value of $SP(n,k)$ is found when $k$ divides $n$ and an upper bound
on $SP(n,k)$ for all $n$ and $k$ is given. These results are restated
below.

\begin{theorem} \label{thm:MMS2}
  Let $k,\ell$ be integers, then 
\[
SP(\ell k,k) =  \binom{\ell k-1}{\ell-1}.
\] 
Moreover, a Sperner partition system meets this bound only if every
class of every partition in the system has exactly $\ell$ elements.
\end{theorem}

For general $n$ only an upper bound on the size of Sperner partition
system was determined in~\cite{MMS}.

\begin{theorem} \label{thm:MMS1} 
  Let $n,k,\ell$ and $r$ be integers with $n=\ell k+r$ and $0 \leq r <
  k$. Then
\[
SP(n,k) \leq \frac{1}{(k-r) + \frac{r(\ell+1)}{n-\ell}} \binom{n}{\ell}.
\]
\end{theorem}

It is not hard to see that the size of Sperner $k$-partitions systems
is non-decreasing with the value of $n$.

\begin{lemma}\label{lem:trivialbound}
For all $n$ and $k$
\[
SP(n-1,k) \leq SP(n,k). 
\]
\end{lemma}
\proof Simply start with a Sperner $k$-partition system on a
$(n-1)$-set with size $SP(n-1,k)$ and add a new $n$-th element to one
of the classes in each partition in the system.\qed

The focus of this paper is to provide exact values
and new bounds for $SP(k,n)$ in cases where either $k$ or $\ell = \lfloor
\frac{n}{k} \rfloor$ is small. Our results are summarized in the list
below:
\begin{enumerate}[(a)]
\item Theorem~\ref{thm:two}: For any integer $\ell$ the exact value of
  $SP(2\ell+1,2)$ is $\binom{2\ell}{\ell-1}$.
\item Theorem~\ref{thm:bound2k+1}: For all $k$, $SP(2k+1,k) \leq 2k$
\item Theorem~\ref{thm:sp2k+1-k}: If $k$ is even, then $SP(2k+1,k) = 2k$
\item Theorems~\ref{thm:2k+2lower} and Theorem~\ref{thm:2k+2upper}:
  If $k\geq 3$, then $2k+1 \leq SP(2k+2,k) \leq 2k+3$.
\item Theorem~\ref{thm:3k-1}:  If $k>2$, then $3k-1 \leq SP(3k-1,k)$.
\item Theorem~\ref{thm:n-k}: For all $n$ and $k$, $SP(n,k) \geq k*SP(n-k,k)$. 
\end{enumerate}

\section{Sperner $2$-partition systems}

In this section we determine the exact size of the largest Sperner
partition system when $k=2$ for all $n$.  If $n$ is even, the
size of the largest Sperner $2$-partition system on an $n$-set is
$\binom{n-1}{\frac{n}{2}-1}$ by Theorem~\ref{thm:MMS2}. If $n = 2 \ell + 1$
Theorem~\ref{thm:MMS1} gives the following bound
\[
SP(n,k) \leq \frac{1}{ 1 + \frac{ (\ell+1)}{\ell+1}}  \binom{2\ell+1}{\ell}
       = \left( 1+ \frac{1}{2\ell} \right) \binom{2\ell}{\ell-1}.
\]
We will show that the exact maximum size of such a Sperner partition
system is $\binom{2\ell}{\ell-1}$.

\begin{theorem} \label{thm:two} If $k=2$ and $n = 2\ell +
  1$, then $SP(n,2) = \binom{n-1}{\ell-1}$. In addition, equality
  holds if and only if each partition has one class of size $\ell$ and
  the other class of size $\ell+1$.
\end{theorem}
\proof It is not hard to construct a Sperner $2$-partition system on
an $n$-set with $\binom{n-1}{\ell-1}$ partitions. Set each partition in
this system to be one set of size $\ell$ that contains $1$ and its
complement (note that the complement will have size $\ell+1$).  It is
clear that this is a Sperner partition system since the classes of the
partitions are all distinct and no class of size $\ell+1$ can contain
one of the classes of size $\ell$ (since all the classes of size $\ell$
contain $1$, but none of the classes of size $\ell+1$ do).

Next we must show that this is the largest size of a Sperner
$2$-partition system. In any such system, any partition must contain
one class with at most $\ell$ elements.  The collection of all such
classes must be an intersecting set system (meaning any two such sets have
non-trivial intersection). If this were not the case, then there
would be two sets $A,B$ in the set system with $A \cap B = 0$, which
implies that $A \subseteq \overline{B}$, which contradicts the
assumption that these sets came from a Sperner partition system.
Therefore, by the Erd\H{o}s-Ko-Rado Theorem~\cite{MR0140419}, there
are at most $\binom{n-1}{\ell-1}$ such sets and $SP(n,2) \leq
\binom{n-1}{\ell-1}$.

Finally, equality holds in the Erd\H{o}s-Ko-Rado Theorem only if each
set in the set system has size exactly $\ell$, which implies a Sperner
$2$-partition system has $\binom{n-1}{\ell-1}$ partitions only if each
partition has one class of size $\ell$ and the other of size $\ell+1$.  
\qed

\section{Sperner $k$-partition systems on a $2k+1$ set}

Next we consider the case where $n = 2k+1$; this is the least value of
$n$ for which the value of $SP(n,k)$ is not known.  Indeed, if $n
=2k$, then the exact value of $SP(n,k)$ is given by
Theorem~\ref{thm:MMS2}. If $n < 2k$, then there must be a class of
size $1$ in any $k$-partition of an $n$-set, thus there cannot be two
Sperner $k$-partitions and $SP(n,k)=1$. In fact, a Sperner partition
system with more than two partitions cannot contain a partition with a
class of size one and, from this fact, it is clear that any Sperner $k$-partition
system on a $2k+1$ set (with more than one partition) must contain only
partitions consisting of $k-1$ classes of size $2$ and one class of
size $3$.

For $n=2k+1$, the bound from Theorem~\ref{thm:MMS1} is
\[
SP(2k+1,k) \leq \frac{1}{(k-1) + \frac{3}{2k-1}} \binom{2k+1}{2}
= \frac{ 4k^3-k} { 2k^2-3k+4 }  \leq 2k+3.
\]
Putting Lemma~\ref{lem:trivialbound} and Theorem~\ref{thm:MMS2}
together produces a lower bound on $SP(2k+1,k)$ of $2k-1$.  The next
theorem gives a lower bound that is the exact value of $SP(2k+1,k)$
for all even values of $k$.

\begin{theorem}\label{thm:bound2k+1} 
 For all integers $k$, it is the case that  $SP(2k+1,k) \leq 2k$.
\end{theorem}
\proof Suppose $\mathcal{P}$ is a Sperner partition system containing
$2k+1$ partitions. Assume that $1$ is the element that appears in the
fewest number of classes of size three.

No element can appear only in classes of size two, since there are
$2k+1$ partitions and $2k+1$ elements in total. Therefore, every
element must appear in at least one class of size three. Since there are
$2k+1$ classes of size three, there must be an element that occurs in no
more than three classes of size three. Thus the element $1$ must belong to one, two
or three classes of size three.

Consider the case where $1$ belongs to exactly three classes of size
three (the cases when $1$ belongs to one or two classes of size three is
similar). The element $1$ must belong to $2k-2$ classes of size two,
which means that there are exactly two elements which do not appear
with $1$ in a class of size two. Thus three classes of size three that
contain $1$ are not distinct.
\qed

Next we will prove that $SP(2k+1,k) \geq 2k$ when $k$ is even by
giving a construction for a Sperner partition system of this size.
This construction is based on a well-known construction for a
1-factorization of the complete graph on an even number of points
(see~\cite[VII.5.5]{MR2246267} for the details of this construction).  We will
start with an example of this construction for $k=8$ and then give the
details of the construction for all even $k$.  

Place $16$ points in a circle and label the points $1$ to $16$; place
a $17$-th point in the center of the circle and label this point by
$\infty$. We will define a single triangle and $k-1$ edges on these
points. This will be the initial partition in the Sperner partition
system. Leaving the labeled points in place, we will rotate the edges
and the triangle to construct a new partition. Rotating this partition
$16$ times will produce $16$ partitions; these will be the partitions
in the Sperner system. We will call this process of rotating an
initial partition to create $2k$ new partitions \textsl{developing}
the initial partition.

For $k=8$ the initial partition is 
\[
\{1,5,9\},\;  \{8,11\},\;\{7, 12\},\;\{6, 13\}, \; \{2,16\},\; \{4, 10\},\; \{3,\infty\},  \;   \{14,15\}
\]
(the method used to construct this initial partition is given in
Theorem~\ref{thm:sp2k+1-k}) this partition is represented in the
following diagram as a triangle and $7$ edges.

\begin{figure*}[h]
\begin{center}
\begin{tikzpicture}[scale=0.6]
\draw  [fill] (0,2.2)  circle (2pt) node[above]{1} 
 --  (0,-2.4) circle (2pt) node [below]{9}  
   -- (2,0) circle (2pt) node [right]{5};
\draw (2,0) --(0,2.2) ;
\draw [fill] (0.7,1.8)  circle (2pt) node[above]{2}
 -- (-0.7,1.8)  circle (2pt) node[above]{16};
\draw [fill] (1.4,1.4)  circle (2pt) node[above right] {3} 
   -- (.25,.25)  circle (2pt) node[above] {$\infty$};
\draw [fill] (1.8,0.8)  circle (2pt) node[right] {4}
-- (-0.7,-2.0)  circle (2pt) node[below left] {10};
\draw [fill] (1.8,-0.8)  circle (2pt) node[right] {6}
   --  (-2,0)  circle (2pt) node[left] {13};
\draw [fill] (1.4,-1.4)  circle (2pt) node[below right] {7} 
-- (-1.8,-0.8)  circle (2pt) node[left] {12};
\draw [fill] (0.7,-2.0)  circle (2pt) node[below right] {8} 
-- (-1.4,-1.4)  circle (2pt) node[left] {11};
\draw [fill]  (-1.8,0.8)  circle (2pt) node[left] {14}
-- (-1.4,1.4)  circle (2pt) node[above left] {15};
\end{tikzpicture}
\end{center}

\caption{The edges and triangle that form the initial partition in the
  Sperner partition system in Figure~\ref{tab:17-8}.\label{fig:initialk8}}
\end{figure*}
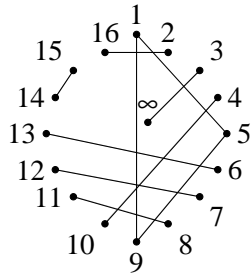

\begin{figure}
\begin{center}
\small{
\begin{tabular}{l|l}
partition & classes\\ \hline
1& $\{1, 5, 9\}, \; \{8, 11\},\;\{7, 12\},\;\{6, 13\},\{2, 16\},\;\{4, 10\},\;\{\infty, 3\},\;\{14, 15\}$ \\ 
2& $\{2, 6, 10\},\;\{9, 12\},\;\{8, 13\},\;\{7, 14\},\;\{3, 1\},\;\{5, 11\},\;\{\infty, 4\},\;\{15, 16\} $ \\ 
3& $\{3, 7, 11\},\;\{10, 13\},\;\{9, 14\},\;\{8, 15\},\;\{4, 2\},\;\{6, 12\},\;\{\infty, 5\},\;\{16, 1\} $ \\ 
4& $\{4, 8, 12\},\;\{11, 14\},\;\{10, 15\},\;\{9, 16\},\;\{5, 3\},\;\{7, 13\},\;\{\infty, 6\},\;\{1, 2\} $ \\ 
5& $\{5, 9, 13\},\;\{12, 15\},\;\{11, 16\},\;\{10, 1\},\;\{6, 4\},\;\{8, 14\},\;\{\infty, 7\},\;\{2, 3\} $ \\ 
6& $\{6, 10, 14\},\;\{13, 16\},\;\{12, 1\},\;\{11, 2\},\;\{7, 5\},\;\{9, 15\},\;\{\infty, 8\},\;\{3, 4\} $ \\ 
7& $ \{7, 11, 15\},\;\{14, 1\},\;\{13, 2\},\;\{12, 3\},\;\{8, 6\},\;\{10, 16\},\;\{\infty, 9\},\;\{4, 5\} $ \\ 
8& $ \{8, 12, 16\},\;\{15, 2\},\;\{14, 3\},\;\{13, 4\},\;\{9, 7\},\;\{11, 1\},\;\{\infty, 10\},\;\{5, 6\} $ \\ 
9& $ \{9, 13, 1\},\;\{16, 3\},\;\{15, 4\},\;\{14, 5\},\;\{10, 8\},\;\{12, 2\},\;\{\infty, 11\},\;\{6, 7\} $ \\ 
10& $ \{10, 14, 2\},\;\{1, 4\},\;\{16, 5\},\;\{15, 6\},\;\{11, 9\},\;\{13, 3\},\;\{\infty, 12\},\;\{7, 8\} $ \\ 
11& $ \{11, 15, 3\},\;\{2, 5\},\;\{1, 6\},\;\{16, 7\},\;\{12, 10\},\;\{14, 4\},\;\{\infty, 13\},\;\{8, 9\} $ \\ 
12& $ \{12, 16, 4\},\;\{3, 6\},\;\{2, 7\},\;\{1, 8\},\;\{13, 11\},\;\{15, 5\},\;\{\infty, 14\},\;\{9, 10\} $ \\ 
13& $ \{13, 1, 5\},\;\{4, 7\},\;\{3, 8\},\;\{2, 9\},\;\{14, 12\},\;\{16, 6\},\;\{\infty, 15\},\;\{10, 11\} $ \\ 
14& $ \{14, 2, 6\},\;\{5, 8\},\;\{4, 9\},\;\{3, 10\},\;\{15, 13\},\;\{1, 7\},\;\{\infty, 16\},\;\{11, 12\} $ \\ 
15& $ \{15, 3, 7\},\;\{6, 9\},\;\{5, 10\},\;\{4, 11\},\;\{16, 14\},\;\{2, 8\},\;\{\infty, 1\},\;\{12, 13\} $ \\ 
16& $ \{16, 4, 8\},\;\{7, 10\},\;\{6, 11\},\;\{5, 12\},\;\{1, 15\},\;\{3, 9\},\;\{\infty, 2\},\;\{13, 14\} $ \\
\end{tabular}
}
\end{center}
\caption{A Sperner $8$-partition system on an $17$-set with $16$ partitions.\label{tab:17-8}}
\end{figure}

Developing this partition yields the Sperner partition system given in
Figure~\ref{tab:17-8}.  The key to this construction lies in the
differences between the elements in the classes of the
partitions. Define the difference between any two points $i$ and $j$
to be simply the minimum of $i-j \mod {2k}$ and $j-i \mod {2k}$ (this
is the shortest distance around the circle between the points), except
if one of $i$ or $j$ is $\infty$, in this case we will say that the
distance is infinity. Addition and subtraction of the vertices will be
modulo $2k$ and this convention will be used in other constructions as well.

The first column of Figure~\ref{tab:17-8} contains the classes of size
three for each partition, note that the differences between the any
two of the elements in any of these classes is either $4$ or $8$. In
the next column the difference between the elements in the class is
always $3$. In the third column the difference is $5$ and in the
fourth column the difference is $7$. Each of the elements in the
classes in the fifth column differ by $2$ and by $6$ in the sixth
column. The pairs in the seventh column all include the center vertex
and the difference between them is $\infty$. The final column contains
only pairs that differ by one. The pairs of points in each column all
have the same differences and no two columns have the same
differences. Thus no class from one column can be a subset of a
partition from another column. By construction, the classes within a
column are all distinct.  Thus this is indeed a Sperner partition
system of the correct size.

The key to this construction is to construct an initial partition for
each value of $k$ so that each of the difference between the elements
in a class of size two are all distinct and none of these differences
are equal to the differences between the elements in the class of size
three. We will say that a partition with this property has the
\textsl{difference property}.  Since there are $k+1$ distinct possible
differences of the vertices in the circle ($1$ to $k$ and $\infty$)
and there are exactly $k-1$ pairs in the partition, the triangle must
only have two distinct differences. Further the difference $k$ must be
contained in the triangle, since if it is covered by an edge, this
edge will be repeated when the partition is developed.

\begin{theorem}\label{thm:sp2k+1-k}
  For $k$ even, there exists a Sperner $k$-partition system on a
  $(2k+1)$-set with $2k$ partitions.
\end{theorem}
\proof Place $2k$ points in a circle and label them $1$ to $2k$, add
another point in the center of the circle and label it $\infty$. We
will give a construction for a set of $k-1$ pairs and one triangle on
these points that will be the first partition in the Sperner
system. We will show that the system constructed by developing this
partition is a Sperner $k$-partition system on a $2k+1$ set with $2k$
partitions.

There are two slightly different constructions, depending on the
value of $k$ modulo four. In both of these cases, we will only give
the construction of the initial partition and show that the partition
has the difference property.

In this construction, triangle is determined and the remaining
differences are covered by edges. The approach is to use a series of
parallel edges to cover the even differences and another series of
parallel edges for the odd differences. To avoid covering the
differences in triangle twice, it is sometimes necessary to break from
this pattern.

\medskip
\noindent{\bf Case 1: $k \equiv 0 \pmod{4}$} \hfill 

The construction of the initial partition in this case is the more
complicated; the partition pictured in
Figure~\ref{fig:initialk8} is based on this construction.

For this partition the triangle is $\{1,k/2+1,k+1\}$, this covers the
differences $k$ and $\frac{k}{2}$. Note that both $k$ and $k/2$ are
even. Next we cover all the odd differences. Use the pair $\{k-i,
k+1+i\}$ to cover the difference $2i+1$ for each $i = 1,...,k/2-1$. To do
this, it is necessary that $k>4$.

It is slightly more complicated to cover the remaining even
differences. First use the pairs $\{1-i,1+i\}$ for $i = 1, ...,
(k-4)/4$ to cover the even differences $2,...,k/2-2$.  The difference
$k/2$ is covered in the triangle. Now we must consider two
subcases. If $k/4$ is even, then place $k+2$ in a class with $k/4+2$
and $\infty$ in a class with $k/4 +1$. This covers the differences
$\infty$ and $3k/4$ (which is even).  If $k/4$ is odd, place $k+2$ in
a class with $k/4+1$ and place $0$ with $k/4+2$.  This will cover the
differences $\infty$ and $3k/4+1$ (which is even). We will use the
pairs $\{1-i,1+i+2\}$ to cover the even differences between $k/2+2$
and $3k/4-2$ or $3k/4-1$ (whichever is even). The even difference
$3k/4$ (or $3k/4 +1$) cannot be covered a second time, so we use the
pair $\{2-3k/8, 1-3k/8\}$ (or $\{2-(3k+4)/8, 1-(3k+4)/8\}$) to skip
this difference --- note that this pair covers the difference
$1$. Finally we use the pairs $\{1-i+2,1+i+2\}$ to cover the remaining
even differences.

We summarize the classes in the initial partition in the following
charts. The first chart is for when $k/4$ even and the second for when it is odd.

\begin{center}
\begin{tabular}{|l|l|l|}\hline
class & difference & values of parameter\\ \hline
$\{1,k/2+1, k+1\}$ & $k$ and $\frac{k}{2}$ & \\
$\{k-i, k+1+i\}$ & $2i+1= 3,5,...,k-1$ & $i = 1,...,k/2-1$ \\
$\{1-i,1+i\}$ &  $2i=2,4,...,k/2-2$ &  $i = 1, ..., (k-4)/4$ \\
$\{k+2, k/4+2\}$ & $3k/4$ & \\
$\{\infty, k/4+1\}$  & $\infty$ & \\
$\{1-i,1+i+2\}$ & $2i=k/2+2,...,3k/4-2$ & $i=k/4,...,3k/8-2$ \\
$\{2-3k/8, 1-3k/8\}$ & $1$ & \\
$\{3-i,i+3\}$ & $2i+2 = 3k/4+2,...,k-2$ & $i=3k/8,...,(k-4)/2$ \\ \hline
\end{tabular}
\end{center}
\medskip

\begin{center}
\begin{tabular}{|l|l|l|}\hline
class & difference & values of parameter\\ \hline
$\{1,k/2+1, k+1\}$ & $k$ and $\frac{k}{2}$ & \\
$\{k-i, k+1+i\}$ & $2i+1 = 3,...,k-1$ & $i = 1,...,k/2-1$ \\
$\{1-i,1+i\}$ &  $2i= 2,...,k/2-2$ &  $i = 1, ..., \frac{k-4}{4}$ \\
$\{k+2, k/4+1\}$ & $3k/4+1$ &\\
$\{\infty, k/4\}$ & $\infty$  &\\
$\{1-i,1+i+2\}$ & $2i+2 = k/2+2,...,3k/4-1$ & $i=k/4,...,\frac{3k+4}{8}-2$ \\  
$\{2-(3k+4)/8, 1-(3k+4)/8\}$ & $1$ & \\
$\{3-i,i+3\}$ & $2i = 3k/4+3,...,k-2$ & $i=\frac{3k+4}{8}+1,...,(k-2)/2$ \\  \hline
\end{tabular}
\end{center}

\medskip
\noindent{\bf Case 2: $k \equiv 2 \pmod{4}$} \hfill 

The triangle in this case is $(1,\frac{k}{2}+1, k+1)$ and it covers
the difference $k$, which is even, and $\frac{k}{2}$, which is odd.  We
will only give a table describing the other classes of the initial
partition.

\begin{center}
\begin{tabular}{|l|l|l|} \hline
class & difference & value of parameter\\ \hline
$\{1,\frac{k}{2}+1, k+1\}$ & $k$ and $\frac{k}{2}$ & \\
$\{k+1-i, k+1+i\}$ & $2i = \{2,4,...,k-2\}$ & $i=1,...,\frac{k-2}{2}$ \\  
$\{\infty,2\}$ & $\infty$ & \\
$\{2+i, 1-i\}$ &  $2i+1=\{3,5, \dots ,\frac{k}{2}-2\}$ & $i=1,...,\frac{k-6}{4}$ \\
$\{\frac{7k+6}{4}, \frac{7k+2}{4}\}$ & $1$ & \\ 
$\{2+i, 1-i-2\}$ &  $2i+3 = \{\frac{k}{2}+2,\frac{k}{2}+4,...,k-1\}$ & $i=\frac{k-2}{4},...,\frac{k-4}{2}$ \\ \hline
\end{tabular}
\end{center}

The above construction does not include $k=4$ or $k=2$. The case of
$k=2$ is included in Theorem~\ref{thm:two}.  Figure~\ref{tab:9-4}
gives an example of a Sperner partition system that shows that
$SP(9,4) =8$. \qed

\begin{figure} \label{fig-9-4}
\begin{center}
\small{
\begin{tabular}{l|l} 
partition & classes \\ \hline
$1$ & $\{1,2\},\{3,9\},\{4,5\},\{6,7,8\}$ \\
$2$ & $\{1,3\},\{2,7\},\{4,6\},\{5,8,9\}$ \\
$3$ & $\{1,4\},\{2,5\},\{3,8\},\{6,7,9\}$ \\
$4$ & $\{1,5\},\{2,8\},\{3,6\},\{4,7,9\}$ \\
$5$ & $\{1,6\},\{2,3\},\{5,7\},\{4,8,9\}$ \\
$6$ & $\{1,7\},\{2,9\},\{3,4\},\{6,5,8\}$ \\
$7$ & $\{1,8\},\{2,4\},\{3,7\},\{5,6,9\}$ \\
$8$ & $\{1,9\},\{2,6\},\{3,5\},\{4,7,8\}$ \\
\end{tabular}
}
\caption{A Sperner $4$-partition system on a $9$-set with $8$ partitions.\label{tab:9-4}}
\end{center}
\end{figure}

\begin{lemma}
  There does not exist a Sperner $3$-partition system on a $7$-set
  with six partitions, in particular, $SP(7,3) = 5$.
\end{lemma}
\proof The lower bound of $SP(7,3) \geq 5$ comes from
Theorem~\ref{thm:MMS2} and Lemma~\ref{lem:trivialbound}.

Assume that there is a Sperner $k$-partition system with six
partitions. Then there are a total of twelve pairs and six triples in
this system.  It is clear that no element can appear only in classes
of size three, for otherwise, by removing this element, we get a Sperner
$3$-partition system on a $6$-set with six partitions.  Therefore, we
know that every element must appear in a class of size two.

There is an average of $3(2k)/(2k+1) = 18/7 <3$ elements which occur
in a triple, so there must be some element that occurs in fewer than
three pairs.  If this element occurs in only two triples, then it occurs
as a pair with four distinct elements which means it cannot occur in
two distinct triples (as there are only two elements that have not
occur with it). Similarly, if it occurs in only one triple then it
occurs as a pair with five distinct elements, so there are not enough
distinct elements to form a triple.

Thus there is an element which occurs only as a pair, never a
triple.  Since there are six partitions, this element must occur with
every other element exactly once.

Assume that there is a second element that never occurs in a triple,
this element would also have to occur with every other element exactly
once as a pair.  Consider the partition that has both of these
elements in a class of size two. There is another class of size two in
this partition, assume that contains the elements $x$ and $y$.  At
most four partitions contain a class of size two that contain either
$x$ or $y$.  Thus there is a partition that contains $x$ and $y$
together in the class of size three. This contradicts the system being
a Sperner partition system.

Thus there is exactly one element that occurs in fewer than three
classes of size three. We can conclude that every other element occurs
in exactly three classes of size three and three classes of size two.

Consider an element $x$ which occurs in exactly three classes of size
three. This element must occur three times as a pair, so there are
only three elements along with $x$ that make up the classes of size
three that contain $x$.  So the classes must look like: $\{x, a, b \}
, \{x,a,c\},\{x,b,c\}$. Further, the elements $a,b$ and $c$ must also
occur in exactly three classes of size three, so using the same
argument, the set $\{a,b,c\}$ is also a triple in the system. This
means that the elements which occur in the classes of size three can
be grouped into collections of four sets. Then the classes are all
triples from these sets of four elements. This means that six must be
divisible by four which is a contradiction.  \qed

\section{Bounds on $SP(2k+2,k)$}

In this section we consider Sperner $k$-partition system on an $n$
where $n = 2k+2$. In this case we give a construction that is similar
to the one used in the previous section to establish a lower bound on
the size of such a Sperner partition system and we have also give an
upper bound for these systems.

\begin{theorem} \label{thm:2k+2lower}
  For all $k>2$, there exists a Sperner $k$-partition system on a
  $(2k+2)$-set with $2k+1$ partitions.
\end{theorem}
\proof Similar to the proof of Theorem~\ref{thm:sp2k+1-k} place the
$2k+1$ points in a circle and label them $1$ to $2n+1$. Place a final
point in the center and label it $\infty$. We will define a of two
triangles and $k$ edges on these points that form a partition that
when developed forms a Sperner partition system.

The initial partition is
\[
\{ 1,2,\infty\}, \{3,2k+1\},\{4,2k\},...,\{k,k+4\},
\{k+1, k+2, k+3\}
\]
This partition is represented in the picture below.

\begin{figure}[h]
\begin{center}
\begin{tikzpicture}[scale=0.6]
\draw [fill](2.5, 5)   circle (2pt) node[above]{1} 
    -- (3.5,5)   circle (2pt) node[above right]{2} 
    -- (3,2.5) circle (2pt) node [right]{$\infty$} -- (2.5,5);
\draw [fill](1.5,4)   circle (2pt) node[above left]{$2k+1$}
   --  (4.5,4)   circle (2pt) node[above right]{3};
\draw [fill](1,3)   circle (2pt) node[left]{$2k$}
   -- (5,3)   circle (2pt) node[right]{4};
\draw [fill] (1,2.5)  node{$\vdots$}; 
\draw [fill] (5,2.5)  node{$\vdots$};

\draw [fill](1.2,1.5)   circle (2pt) node[left]{$k+4$}
   -- (4.8,1.5)   circle (2pt) node[right]{$k$};

\draw [fill](4.2,.5)   circle (2pt) node[right]{$k+1$}
 --  (3,0)   circle (2pt) node[below]{$k+2$}
 -- (1.8,.5)   circle (2pt) node[left]{$k+3$} -- (4.2,.5);
\end{tikzpicture}
\end{center}
\caption{The edges and triangle that form the initial partition for a
  Sperner $4$-partition system on an $11$-set}      
\end{figure}
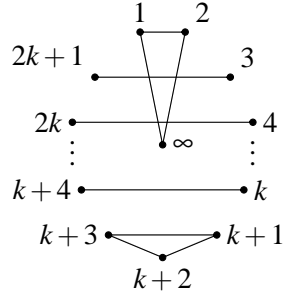

The triangles use only the differences $1$, $2$ and $\infty$ and the
pairs uses difference $3$,$4$,...,$k$ (this first $\lceil \frac{k}{2}
\rceil$ pairs cover the odd differences and the last pairs cover the
even differences). Since the differences in the edges are not contained
in the differences in the triangles and the triangles contain a
distinct set of differences, developing this partition produces a
Sperner partition system.\qed

\begin{theorem} \label{thm:2k+2upper}
  For all $k \geq 3$, the value of $SP(2k+2,k)$ is no more than
  $2k+3$.
\end{theorem}
\proof 
All classes in a Sperner $k$-partition of a $(2k+2)$-set must have
size $2$, $3$ or $4$. We will call classes of size $3$ or $4$ \textsl{large
  classes} and classes of size $2$ \textsl{small classes}.

Let $n=2k+2$ and suppose there is a Sperner $k$-partition system on an
$n$-set that has $n+1=2k+4$ partitions.  Each of the $2k+4$ partitions
must have at least $k-2$ small classes. The average number of times an
element occurs in a small class is $2(2k+4)(k-2)/(2k+2) =
(2k^2-8)/(k+1)$ which is greater than $2k-4$ as $k > 2$. Therefore,
some element, say element $1$, occurs with frequency at least $2k-3$
in the small classes. (Note that $1$ cannot occur in more that $2k+1$
classes of size one.)

First, assume that this frequency is exactly $2k-3$. In the remaining
seven partitions the element $1$ must be in a large class.  These
large classes must contain only elements that do not occur with
element $1$ in the small classes, since no class can be contain in
another.  There are only $(2k+2)-(2k-2) = 4$ such elements, so at most
six large classes containing element $1$ can be formed with these $4$
elements, contradicting the fact we need seven such large classes. A
similar argument works if the element $1$ occurs in more than $2k-3$
small classes.
\qed

The bound in Theorem~\ref{thm:2k+2upper} is a better upper bound than
the bound that can be found from Theorem~\ref{thm:MMS1}.

\section{A lower bound on $SP(3k-1,k)$}

In this section, we will give a construction for a Sperner
$k$-partition system on an $n$-set where $n = 3k-1$.  This
construction produces a system in which all the partitions have $k-1$
classes of size $3$ and one of size $2$. This construction is similar
to the ones given in the proofs of Theorem~\ref{thm:sp2k+1-k} and
Theorem~\ref{thm:2k+2lower}.  Place $3k-1$ points in a circle
and label them $1$ to $3k-1$. We will define $k-1$ distinct triangles
on these points and one pair --- this will be the first partition in
our system. Developing this partition will produce the other
partitions in the Sperner partition system.  As before, we will start
with an example of the construction and then give the details.

\begin{figure}[h]
\begin{center}
\begin{tikzpicture}[scale=0.6]
\draw [fill](4.2,1)   circle (2pt) node[below]{6} ;
\draw [fill](1,2)   circle (2pt) node[left]{9} 
   -- (5,2)  circle (2pt) node[below right]{5} --
      (3,0.5)   circle (2pt) node[below]{7} -- (1,2);
\draw [fill](1,3)   circle (2pt) node[left]{10}
   -- (5,3)   circle (2pt) node[right]{4}
   -- (1.8,1)   circle (2pt) node[left]{8} -- (1,3);
\draw [fill](1.5,4)   circle (2pt) node[above left]{11}
   --  (4.5,4)   circle (2pt) node[above right]{3}
    -- (4.2,1)   circle (2pt) node[below]{6}  -- (1.5,4) ;
\draw [fill](2.5, 5)   circle (2pt) node[above]{1} 
    -- (3.5,5)   circle (2pt) node[above right]{2} ;
\end{tikzpicture}
\end{center}
\caption{The edges and triangle that form the initial partition for a
  Sperner $4$-partition system on an $11$-set}      
\end{figure}

Note that the $(3,6,11)$ triangle covers differences $3$ and $5$, the
$(4,8,10)$ triangle covers the differences $5$, $4$ and $2$ and the
$(5,7,9)$ triangle covers the differences $4$ and $2$. Since each
triangle covers a different set of differences, rotating one triangle
will never produce one of the other triangles. Further, none of the
triangles includes an edge of length one, so the edge will never occur
in one of the triangles. With this construction, the first partition
is $\{1,2\},\{3,6,11\},\{4,8,10\},\{5,7,9\} $ and the rest of the
system is given in Figure~\ref{tab:11-4}.

\begin{figure}[h]
\begin{center}
\small{
\begin{tabular}{l|l} 
partition & classes \\ \hline
1 & $\{1,2\},\;\{3,6,11\},\;\{4,8,10\},\;\{5,7,9\} $ \\
2 & $\{2,3\},\;\{4,7,1\},\;\{5,9,11\},\;\{6,8,10\} $ \\
3 & $\{3,4\},\;\{5,8,2\},\;\{6,10,1\},\;\{7,9,11\} $ \\
4&$\{4,5\},\;\{6,9,3\},\;\{7,11,2\},\;\{8,10,1\} $ \\
5&$\{5,6\},\;\{7,10,4\},\;\{8,1,3\},\;\{9,11,2\} $ \\
6&$\{6,7\},\;\{8,11,5\},\;\{9,2,4\},\;\{10,1,3\} $ \\
7&$\{7,8\},\;\{9,1,6\},\;\{10,3,5\},\;\{11,2,4\} $ \\
8&$\{8,9\},\;\{10,2,7\},\;\{11,4,6\},\;\{1,3,5\} $ \\
9&$\{9,10\},\;\{11,3,8\},\;\{1,5,7\},\;\{2,4,6\} $ \\
10&$\{10,11\},\;\{1,4,9\},\;\{2,6,8\},\;\{3,5,7\} $ \\
11&$\{11,1\},\;\{2,5,10\},\;\{3,7,9\},\;\{4,6,8\} $ \\
\end{tabular}
}
\end{center}
\caption{A Sperner $4$-partition system on an $11$-set with $11$ partitions.\label{tab:11-4}}
\end{figure}

As in Theorem~\ref{thm:sp2k+1-k}, we will give the initial partition
for every $k$ and prove that developing this partition produces a
Sperner $k$-partition system on a $(3k-1)$-set.

\begin{theorem}\label{thm:3k-1}
For any integer $k > 3$
\[
SP(3k-1,k) \geq 3k-1.
\] 
\end{theorem}
\proof 
Start by defining $k$ pairs, 
\[
\{1,3k-1\}, \; \{2,3k-2\}, \; \{3,3k-3\}, ... ,\{k,2k\}.
\] 
The first pair, $\{1, 3k-1\}$ will be the class of size two in the
partition, to each of the other pairs we add one more vertex. Add the
vertex $k+1$ to the edge $\{2,3k-2\}$. To the edge $\{i,3k-i\}$ add
the vertex $2k+2-i$ where $i = 3,...,k-2$.  Let the initial partition
be the edge $\{1, 3k-1\}$, the triangle $\{2,k+1, 3k-2\}$ and the
$k-2$ triangles of the form $\{i,2k+2-i, 3k-i\}$ where $i =
3,...,k-2$.
 
We claim that if this partitions is developed, it will produce a
Sperner partition system.  The triangle $\{2,k+1, 3k-2\}$ covers the
differences $k-1$, $3$ and $k+2$ and each triangle of the form
$\{i,2k+2-i, 3k-i\}$ covers the differences $2k+2-2i$, $2i-1$ and
$k-2$. Note that the value of $2k+2-2i$ is always even and the value
of $2i-1$ is always odd. This means that rotating one of these
triangles $3k-1$ times will produce $3k-1$ distinct new triangles.
Further, if $k>3$, then none of the triangles includes a
difference of $1$, neither the edge $\{1,3k-1\}$, nor any edge
developed from it, are contained in one of the triangles or in any of
the triangles formed when the partition is developed.

We need to prove that no class of size three is repeated when the
initial partition is developed. Assume that two triangles formed by
rotating the distinct triangles $\{i,2k+2-i, 3k-i\}$ and $\{j,2k+2-j,
3k-j\}$ for some values of $i$ and $j$ are the same. Then $\{2k+2-2i,
2i-1, k-2\}=\{2k+2-2j, 2j-1, k-2\}$ since the triangles would have to
cover the same differences. Since $2k+2-2i$ is even and $2j-1$ is odd,
this means that $i=j$. Thus we can conclude that when the first
partition is developed, no triangles of the form $(i,2k+2-i, 3k-i)$
are repeated.

Finally, we only need to show that rotating any triangle of the form
$\{i,2k+2-i, 3k-i\}$ will not produce the triangle formed by rotating
$\{2,k+1, 3k-2\}$. Indeed, if this were to happen, then the differences
formed by the two triangles would have to be the same, namely
$\{2k+2-2i, 2i-1, k-2\}=\{3,k-1,k+2\}$ for some $i$. For $k>1$ it is
not possible for $k-2$ to be equal to either $k-1$ or $k+2$ modulo
$3k-1$. Thus $k-2 =3$ and $k=5$. In this case there will be two
triangles that both form the differences $\{3,4,7\}$ but the
differences will occur in a different order (meaning that the
triangles will have different orientations) so developing these
triangles will not produce repeated triangles.  \qed

The above construction does not give an example for $k=3$ (the
case where $k=2$ is included in Theorem~\ref{thm:two}).
Figure~\ref{tab:8-3} gives a Sperner $3$-partition system on a set of size $8$.

\begin{figure} \label{fig-8-3}
\begin{center}
\small{
\begin{tabular}{l|l} 
partition & classes \\ \hline
$1$ & $\{0,1\},\{2,3,4\},\{5,6,7\}$ \\
$2$ & $\{0,2\},\{1,3,5\},\{4,6,7\}$ \\
$3$ & $\{0,3\},\{1,4,6\},\{2,5,7\}$ \\
$4$ & $\{1,7\},\{0,4,6\},\{2,3,5\}$ \\
$5$ & $\{3,7\},\{0,5,6\},\{1,2,4\}$ \\
$6$ & $\{2,6\},\{0,5,7\},\{1,3,4\}$ \\
$7$ & $\{3,6\},\{0,4,7\},\{1,2,5\}$ \\
$8$ & $\{4,5\},\{0,6,7\},\{1,2,3\}$ \\
\end{tabular}
}
\caption{A Sperner $3$-partition system on an $8$-set with $8$ partitions.\label{tab:8-3}}
\end{center}
\end{figure}

It is not known how close to optimal the Sperner partition system constructed in
Theorem~\ref{thm:3k-1} is. The upper bound on $SP(3k-1,k)$ from
Theorem~\ref{thm:MMS1} is $\frac{1}{2}\binom{3k-1}{k}$; this is much
larger than the size of the system that we can construct, but we
suspect that the bound from Theorem~\ref{thm:MMS1} is not close to the
optimal size in this case.

\section{Further Work}

Determining the size of Sperner $k$-partition systems on an $n$-set is
difficult when $k$ does not divide $n$ and there is clearly still much work to
be done. We will end with two notes about general Sperner partition systems.

First, lemma~\ref{lem:trivialbound} gives a trivial way to construct a larger
Sperner partition system from a smaller one. We offer a better way to
do this.

\begin{theorem} \label{thm:n-k}
For all integers $n$ and $k$ it is the case that $SP(n,k) \geq k*SP(n-k,k)$.
\end{theorem}

\proof There are exactly $k$ permutations of the $k$ element set
$\{n-k+1,n-k+2,...,n\}$ such that no two of them have a common element
in the same position. (To see this, simply consider the rows of a
Latin square of order $k$.) Take a Sperner $k$-partition system on an
$(n-k)$-set and order the classes of each partition. For each
partition create $k$ new partitions by adding the $i$-th element in
each permutation to the $i$-th class. It is easy to see that these
$SP(n-k,k) \times k$ partitions is a Sperner partition with parameters
$n,k$.  \qed

Our second note has to do with the structure of Sperner partition
systems.  A $k$-partition of an $n$-set is called \textsl{almost
  uniform} if every class of the partition has size $\lfloor
\frac{n}{k} \rfloor$ or $\lceil \frac{n}{k} \rceil$. A $k$-partition
system is an \textsl{almost-uniform partition system} if every
partition in the system is almost uniform. In~\cite{MMS} it was
conjectured that the largest Sperner $k$-partition system on an
$n$-set is an almost-uniform partition system. Theorem~\ref{thm:two}
confirms this conjecture for the case where $k=2$.  This conjecture is
not true in general, there are examples of maximum Sperner partition
systems that are not almost uniform. In Figure~\ref{fig-10-4}) an
example of Sperner $4$-partition system on a $10$-set is given, a
complicated counting argument together with a computer search revealed
that this is the largest such system.

\begin{figure} \label{fig-10-4}
\begin{center}
\small{
\begin{tabular}{l|l} 
partition & classes \\ \hline
$1$ & $\{0,1\}, \{2,3\}, \{4,5\}, \{6,7,8,9\} $ \\
$2$ & $\{0,2\}, \{1,3\}, \{4,6,7\}, \{5,8,9\}$ \\
$3$ & $\{0,4\}, \{1,5\}, \{2,6,7\}, \{3,8,9\} $ \\
$4$ & $\{0,6\}, \{4,8\}, \{1,2,7\}, \{3,5,9\} $ \\
$5$ & $\{0,7\}, \{4,9\}, \{1,2,6\}, \{3,5,8\}$ \\
$6$ & $\{1,9\}, \{5,6\}, \{0,3,8\}, \{2,4,7\}$ \\
$7$ & $\{2,9\}, \{3,6\}, \{0,5,8\}, \{1,4,7\}$ \\
$8$ & $\{2,8\}, \{5,7\}, \{0,3,9\}, \{1,4,6\}$ \\
$9$ & $\{1,8\}, \{3,7\}, \{0,5,9\}, \{2,4,6\}$ \\
$10$& $\{2,5\}, \{3,4\}, \{0,8,9\}, \{1,6,7\}$ \\
\end{tabular}
}
\caption{A $(10,4)-SP$ with $10$ partitions.}
\end{center}
\end{figure}

Perhaps a better conjecture would be that the maximum size of a
Sperner partition system can always be achieved by an almost-uniform
system.

%%%%%%%%%%%%%%%%%%%%

\end{document}